\documentclass[fleqn]{mat01}
\usepackage{times,mathtimy,amssymb,latexsym}
\begin{document}

\setcounter{page}{259}
\firstpage{259}

\def\d{\mbox{\rm d}}

\newtheorem{theore}{Theorem}
\renewcommand\thetheore{\arabic{section}.\arabic{theore}}
\newtheorem{theor}[theore]{\bf Theorem}
\newtheorem{rem}[theore]{Remark}
\newtheorem{lem}[theore]{Lemma}
\newtheorem{note}[theore]{Note}
\newtheorem{coro}[theore]{\rm COROLLARY}
\newtheorem{propo}[theore]{\rm PROPOSITION}
\newtheorem{definit}[theore]{\rm DEFINITION}
\newtheorem{exampl}[theore]{Example}
\newtheorem{case}{Case}

\renewcommand\thecase{\it {\rm {\it \roman{case}}}}

\def\corol{\trivlist \item[\hskip \labelsep{COROLLARY.}]}
\def\noteproof{\trivlist \item[\hskip \labelsep{\it Note added in Proof.}]}

\def\defini{\trivlist\item[\hskip\labelsep{DEFINITION}]}
\def\exammp{\trivlist\item[\hskip\labelsep{\it Examples.}]}

\newtheorem{step}{Step}

\renewcommand\theequation{\arabic{equation}}

\def\a{\alpha }
\def\g{\gamma }
\def\s{\sigma }
\def\p{\parallel }
\def\. {\cdot}
\def\va{\varphi}
\def\G{\Gamma}
\def\n{\nabla }
\def\d{\partial }
\def\l{\left }
\def\r{\right }
\def\z{\zeta}
\def\d{\partial }
\def\sx{{\mathrel{\ni\!\in}}}

\title{Higher order Hessian structures on manifolds}

\markboth{R~David Kumar}{Higher order Hessian structures on manifolds}

\author{R~DAVID KUMAR}

\address{Department of Collegiate Education, Government of Andhra Pradesh,
Hyderabad~500~001, India\\
\noindent E-mail: rdkumar1729@yahoo.com}

\volume{115}

\mon{August}

\parts{3}

\pubyear{2005}

\Date{MS received 9 February 2005; revised 9 June 2005}

\begin{abstract}
In this paper we define $n$th order Hessian structures on
manifolds and study them. In particular, when $n = 3$, we make a
detailed study and establish a one-to-one correspondence between
{\it third-order Hessian structures} and a {\it certain class of
connections} on the second-order tangent bundle of a manifold.
Further, we show that a connection on the tangent bundle of a
manifold induces a connection on the second-order tangent bundle.
Also we define second-order geodesics of special second-order
connection which gives a geometric characterization of symmetric
third-order Hessian structures.
\end{abstract}


\maketitle

\section{Introduction}

Let $M$ stand for a Hausdorff, paracompact, smooth $(C^{\infty})$
manifold modelled on a finite dimensional Banach space {\bf E}. Let
$T_{m}M$ denote the tangent space to $M$ at $m$ \cite{2}. We recall Lang's
characterization of $T_{m}M$ to fix our notation \cite{10}. For each chart
$(U, \varphi)$ at $m$, consider triples of the form $(U,\varphi,{\bf
e})$ with ${\bf e} \in {\bf E}$. If $(U, \varphi)$ and $(V, \psi)$ are
two charts at $m$, define $(U, \varphi, {\bf e})$ and $(V, \psi,
\bar{\bf e})$ to be equivalent if $D(\psi \circ \varphi^{-1})(\varphi m)
\cdot {\bf e} = \bar{\bf e}$. An equivalence class {\bf v} is called a
tangent vector at $m$ and the set of all tangent vectors at $m$ is
denoted by $T_{m}M$. If ${\bf v} \in T_{m}M$ and $(U,\varphi,{\bf e})
\in {\bf v}$, we write ${\bf v}_{\varphi} = {\bf e}$ and $T\varphi ({\bf
v})= (\varphi m, {\bf e})$. If ${\bf v}$, ${\bf w} \in T_{m}M$ and $a,b$
are real numbers, we define $a{\bf v} + b{\bf w} = T\varphi^{-1}
(\varphi m, a{\bf v}_{\varphi} + b{\bf w}_{\varphi})$. This defines a
linear structure on $T_{m}M$ isomorphic to {\bf E} and is independent of
charts. Let $TM = \bigcup_{m \in M} T_{m} M$ and let $TU = \bigcup_{m
\in U} T_{m}M$. Given a chart $(U,\varphi)$ at $m, T \varphi\hbox{\rm :}\ TU
\rightarrow \varphi (U)\times {\bf E}$ is a bijection and $\{(TU,
T\varphi )\hbox{\rm :}\ (U,\varphi)~~ \mbox{is a chart on}~~ M \}$ defines a vector
bundle structure on $TM$. Finally, if $M = U'$ is an open subset of
${\bf E}$, then $TU'$ can be identified with $U'\times {\bf E}$.

Let ${\cal F}(M)$ denote the space of smooth real-valued functions
on $M$. More generally, for any open subset $U$ of $M$, ${\cal
F}(U)$ will stand for the space of smooth real-valued functions on
$U$. If $(U,\varphi)$ is a chart on $M$ and $f \in {\cal F}(U)$ we
write $f_{\varphi} = f \circ \varphi^{-1}$. If $m \in U$, we
define $df(m)\in T^{*}_{m}M$ by $df(m)({\bf v}) =
Df_{\varphi}(\varphi m)\cdot {\bf v}_{\varphi}$ where ${\bf
v}_{\varphi}$ is a local representation of {\bf v} in the chart
$(U, \varphi)$. This definition is independent of charts. For each
fixed {\bf v} in $T_{m}M$, the map $f \mapsto df(m)({\bf v})$ is a
derivation in $f$ and every derivation is of this form for a
unique $\bf v$. Suppose $X\hbox{\rm :}\ m \mapsto X(m)$ is a cross section of
$TM$ over $M$. Given a  chart $(U, \varphi)$, write
$X_{\varphi}(\varphi m) = [X(m)]_{\varphi}$  for $m \in U$. Thus
$T\varphi (X(m)) = (\varphi m, X_{\varphi}(\varphi m))$ and
$X_{\varphi}\hbox{\rm :}\ \varphi U \rightarrow {\bf E}$. We call $X$ a
(smooth) vector field on $M$ if $X_{\varphi}$ is smooth for all
charts $(U, \varphi)$. We denote the set of all vector fields on
$M$ by $ \sx(M) $. If $f \in {\cal F}(M)$ and $X \in \sx(M)$, we
define $X(f)(m) = df(m) \cdot X(m)$. Then $X(f)\in {\cal F}(M)$
and the map $X$ is a derivation on ${\cal F}(M)$. Conversely, every
derivation of ${\cal F}(M)$ is induced in this way by a unique $ X
\in  \sx(M)$. Suppose $M = U'$ is open in ${\bf E}$. Then vector
fields $X$ and $Y$  on $U'$  can be identified with smooth maps
from $U'$ to $\bf E$. We define $DY \cdot X$ to be the smooth map
given by $(DY \cdot X)(u) = DY(u) \cdot X(u)$.

We can see that a (second-order) Hessian structure on a manifold
is equivalent to a symmetric connection on a manifold. In \cite{3},
this equivalence is established by first identifying Hessian
structures with sprays and then sprays with symmetric connections.
The difficulty of their proof is in the choice of defining a
connection in terms of connection forms on the bundle of bases for
the tangent spaces. When we take a connection $\n _{X}Y$ to be
given by Koszul's definition, we see that Hessian structures and
symmetric connections can be directly related to each other.

Before proceeding further, we state certain results relating to
higher order derivatives on manifolds. For $m \in M$, let
$F_{m}^{+}$ denote the space of all those real-valued smooth
functions whose domain is an open subset of $M$ containing the
point $m$, $F_{m}$ stands for the elements of $F_{m}^{+}$ which
vanish at $m$ and $F_{m}^{n}$ for the space of all sums of
products of $n$ elements of $F_{m}$. $f \in F_{m}^{n}$ if and only
if $f_{\va}$ has vanishing partial derivatives of order up to and
including $(n-1)$ at $m$, for some and hence any chart $(U,\va )$.
If $f \in F_{m}^{n}$ then we say $f$ is $(n-1)$-flat at $m$.

\begin{propo}$\left.\right.$\vspace{.5pc}

\noindent Let $f \in F_{m}^{n}$ and let $X_{1}, \ldots, X_{n} \in \sx(M)$. Then
\begin{enumerate}
\renewcommand\labelenumi{\rm (\alph{enumi})}
\leftskip .15pc
\item $[X_{1}\cdot X_{2}\ldots X_{n} \cdot f] (m)$ is defined
independently of the order of $X_{i}${\rm '}s.

\item If ${\bf v}_{1}, \ldots , {\bf v}_{n} \in T_{m}M$ and
$X_{1}, \ldots, X_{n} \in \sx(M)$ such that $X_{i}(m) = {\bf v}_{i}${\rm ,} then the $n$-th
derivative of $f$ at $m${\rm ,}
\begin{equation*}
D^{n}f(m)({\bf v}_{1}, \ldots, {\bf v}_{n}) \ \substack{{\rm def}\\ {=}} \
[X_{1}\cdot X_{2} \ldots X_{n}\cdot f](m)
\end{equation*}
is well-defined and is symmetric in the ${\bf v}_{i}${\rm '}s.

\item If $(U, \va)$ is a chart at $m${\rm ,} then
\begin{equation*}
D^{n}f(m)({\bf v}_{1}, \ldots, {\bf v}_{n}) = D^{n}f_{\va}(\va m)({\bf v}_{1\va},
\ldots, {\bf v}_{n\va })(\va m),
\end{equation*}
where $D^{n}f_{\va}$ is the $n$-th order Frechet derivative of $f_{\va}$ at $\va m$.\vspace{-.4pc}
\end{enumerate}
\end{propo}

If $f\in F_{m}^{n}$, then the $n$th order derivative of $f$ is
well-defined as a symmetric $n$-multilinear mapping on $T_{m}M$. We
define a $n$th order Hessian structure $H^{n}$ (not necessarily
symmetric), for every $n \in {\bf N}$ as a map $H^{n}\hbox{\rm :}\
f\mapsto H^{n}f$ which is real linear in $f$ and associates with every
$f\in {\cal F}(M)$, an $n$th order covariant tensor $H^{n}f$ on $M$ such
that if $f\in F_{m}^{n}$, $H^{n}f(m)$ is the $n$th order derivative of
$f$ at $m$. If $\n $ is a connection, then the definition of $\n ^{n}f$
given by $\n f = df$, $\n ^{n}f = \n (\n ^{n-1}f)$ where the action of
$\n $ denotes the covariant derivative acting on covariant tensors
defines a Hessian structure $H^{n}$, for all $n$. Further we show that
$\n ^{2}f$ and $\n ^{3}f$ are symmetric (for all $f$) if and only if $\n
$ has torsion and curvature zero. In such a case $\n ^{n}f$ is
symmetric, for all $n$ and for all $f$. Thus we show how to define
higher-order derivatives of functions on manifolds with a connection. It
should be interesting to relate these higher derivatives to differential
operator on the manifold.

Ambrose, Palais and Singer in \cite{3} proved that, given any $H (=
H^{2})$ there exists a connection $\n $ such that $H^{2} = \n ^{2}$. We
now raise the following question. Given $H^{3}$, does there exist a
connection $\n $ such that $H^{3} = \n ^{3}$? This is not true in
general. We prove that, just as a second-order Hessian structure arises
from a connection on the tangent bundle, every third-order Hessian
structure on a manifold arises from a connection on the {\it
second-order tangent bundle}. We characterize those $H^{3}$ which arise
as $\n ^{3}$. In this process, we show that a connection $\n $ on the
tangent bundle of a manifold induces a connection $\widetilde{\n }$ on
the second-order tangent bundle. As far as we are aware, this
observation is new. Also we introduce and discuss second-order geodesics
which are related to third-order Hessian structures. We show that if
$\widetilde{\n}$ is a connection on the second-order tangent bundle
induced from a connection $\n $ on the tangent bundle of a manifold,
then every (first-order) geodesic of $\n $ is a second-order geodesic of
$\widetilde{\n }$.

\subsection{\it Second-order tangent bundle}

We first introduce the concept of a second-order tangent vector at a
point $m$ of the manifold $M$. We follow the notation of \cite{3}. Let
$F_{m}^{+}$ denote the space of all those real-valued smooth functions
whose domain is an open subset of $M$ containing the point $m$. Let
$F_{m}$ stand for the elements of $F_{m}^{+}$ which vanish at $m$ and
$F_{m}^{k}$ for the space of all sums of products of $k$ elements of
$F_{m}$. A {\it tangent vector} to $M$ at $m$ can be considered to be a
linear functional on $F_{m}^{+}$ which vanishes on constant functions
and on $F_{m}^{2}$ ( i.e., a tangent vector at $m$ kills on 1-flat
functions at $m$). We define a {\it second-order tangent vector} $t$ to
$M$ as a linear form on $F_{m}^{+}$ which vanishes on constant functions
and on $F_{m}^{3}$ (i.e., $t$ kills on 2-flat functions at $m$). The
second-order tangent vectors form a linear subspace which we denote by
$T_{m}^{(2)}M$. Observe that $T_{m}M \subset T_{m}^{(2)}M$. If $(U,x)$
is a chart at $m$, it is not difficult to show that the functionals
$\frac {\partial }{\partial x^{i}}(m)$ and $\frac{\partial^{2}
}{\partial x^{i}\partial x^{j}}(m) (i\leq j )$ constitute a basis for
$T_{m}^{(2)}M$. Further, any $t\in T_{m}^{(2)}M$ can be written uniquely
in the form $t=\Sigma_{i} t^{i} \frac {\partial }{\partial x^{i}}(m) +
\Sigma_{i\leq j}t^{ij}\frac{\partial^{2} }{\partial x^{i}\partial
x^{j}}(m)$ with the condition that $t^{ij}=t^{ji}$. The span of the
vectors $\frac {\partial }{\partial x^{i}}(m)$ is $T_{m}M$. Let the
vectors $\frac{\partial^{2} }{\partial x^{i}\partial x^{j}}(m)$ span the
subspace $(T_{m}M)^{c}$. Clearly then $T_{m}^{(2)}M = T_{m}M \oplus
(T_{m}M)^{c}$. Moreover $(T_{m}M)^{c}$ is not defined independently of
the chart and there is no canonical way of choosing a complement of
$T_{m}M$ of $T_{m}^{(2)}M$. This leads to the concept of a `dissection'
of second-order tangent bundle.

\begin{defini}$\left.\right.$\vspace{.5pc}

\noindent {\rm A smooth second-order tangent vector field on $M$ is a
map $\xi$ which associates to each point $m$ of $M$ a second-order
tangent vector $\xi_{m}\in T_{m}^{(2)}M$ such that for every $f \in
{\cal F}(M)$, the function $\xi f $ defined by $\xi( f)(m)=\xi_{m}(f)$
is a smooth function on $M$.}
\end{defini}

\begin{exammp}
{\rm $\frac {\partial }{\partial x^{i}}$ and $\frac{\partial^{2}
}{\partial x^{i}\partial x^{j}}$ are smooth second-order vector fields
on their domains.}\vspace{.5pc}
\end{exammp}

Actually, we consider the problem of extending the theory of
second-order structures to Banach manifolds. On Banach manifolds
global chart-free methods are generally not available. One has to
choose charts taking values in open sets of Banach spaces and use
Frechet calculus methods. For this purpose, we have chosen
coordinate-free methods, but use charts combined with Frechet
calculus on Banach spaces. This method works both in finite and
infinite dimensional cases and the proofs are exactly the same as
far as the computational details are concerned. It is for
establishing smoothness and taking into account such infinite
dimensional phenomena as non-reflexivity etc, that one has to
work harder in the case of Banach manifolds. So we use the
language of Frechet calculus, but in a finite dimensional setting.
For this reason our calculations become little complicated and not
easily readable. Our definition of second-order tangent bundle
seems to be new (it is a coordinate approach). Now we define the
second-order tangent bundle of a manifold in the following way.

If $M = {\bf R}^{2}$, the second-order tangent vectors at a point are
the linear functionals of the form $\displaystyle \Sigma a_{i}
(\partial/\partial x_{i}) + \Sigma b_{ij} (\partial^{2}/\partial
x_{i}\partial x_{j})$. Calculating how these functionals transform under
a change of coordinates and rewriting the transformation law in a
coordinate-free language, we are led to the definition of second-order
tangent vector at a point of a manifold.

Let ${\bf E}\Delta {\bf E}$ denote the symmetric tensor product of ${\bf
E}$ with itself. We shall identify symmetric bilinear maps from ${\bf
E}\times {\bf E}$ to ${\bf E}$ with linear maps from ${\bf E}\Delta {\bf
E}$ to ${\bf E}$. If $A$ and $B$ are (linear) endomorphisms of ${\bf
E}$, then $A \Delta B$ denotes the endomorphism of ${\bf E} \Delta {\bf
E}$ satisfying $(A\Delta B)({\bf e}_{1} \Delta {\bf e }_{2}) = A {\bf
e}_{1} \Delta B {\bf e}_{2}$. Let ${\bf E}^{(2)}$ stand for
$({\bf E} \Delta {\bf E}) \oplus {\bf E}$ and we denote a typical element
of it by ${\bf x\oplus e}$.

Let $m \in M$. Consider a chart $(U,\varphi)$ at $m$ and triples of the
form $(U, \varphi, {\bf x \oplus e})$. Two such triples $(U, \varphi,
{\bf x \oplus e})$ and $(V, \psi , {\bf y \oplus f})$ are said to be
equivalent if
\begin{equation*}
{\bf y} = [DF(u) \Delta  DF(u)] \cdot {\bf x}
\end{equation*}
and
\begin{equation*}
{\bf f} = DF(u)\cdot {\bf e} + D^{2} F(u) \cdot {\bf x},
\end{equation*}
where $F = \psi \circ \varphi^{-1}$ and $u = \varphi m$. It can be
easily checked that we have really defined an equivalence relation and
that equivalence preserves linearity. Equivalence classes are called
{\it second-order tangent vectors} at $m$. The set of all such equivalence classes will be
denoted by $T^{(2)}_{m}M$. A typical element of $T^{(2)}_{m}M$ will be
denoted by {\bf t}.

If $(U, \varphi, {\bf x \oplus e}) \in {\bf t}$, we write ${\bf
t}_{\varphi}= {\bf x \oplus e}$ and $T^{(2)}\varphi ({\bf t}) = (\varphi
m,{\bf t}_{\varphi})$. As usual $T^{(2)}M =\bigcup_{m\in
M}T^{(2)}_{m}M$, $T^{(2)}U =\bigcup_{m\in U}T^{(2)}_{m}M$ and
$T^{(2)}\varphi\hbox{\rm :}\ T^{(2)}U \rightarrow \varphi U \times {\bf
E}^{(2)}$ is a bijection. The collection $\{(T^{(2)}U,T^{(2)}\varphi
)\hbox{\rm :}\ (U,\varphi )$ is a chart on $M \}$ defines a vector
bundle structure on $T^{(2)}M$. If $M = U'$ is an open subset of ${\bf
E}$ then $T^{(2)}U'$ can be identified with $U'\times {\bf E}^{(2)}$.

For $m \in M$, let $(T^{(2)}_{m}M)^{*}$ denote the dual of
$T^{(2)}_{m}M$ and let $(T^{(2)}M)^{*}$ denote the dual bundle.

Suppose $f$ is a smooth function defined in a neighborhood of $m$, ${\bf
t} \in T^{(2)}_{m}M$, $(U,\varphi)$ is a chart at $m$, ${\bf t}_{\varphi
} = {\bf x \oplus e}$ and $\varphi m = u$. We define the second-order
differential of $f$ as $d^{(2)}f(m)({\bf t}) \ \substack{{\rm def}\\ {=}} \ {\bf
t}(f)(m) \ \substack{{\rm def}\\ {=}} \ D^{2}f_{\varphi }(u)\cdot {\bf x} +
Df_{\varphi }(u)\cdot {\bf e}.$ It is easy to check that
$d^{(2)}f(m)({\bf t})$ is well-defined and that $d^{(2)}f(m) \in
(T^{(2)}_{m}M)^{*}$. Moreover if $f \in {\cal F}(U)$, $d^{(2)}f$ is a
smooth cross-section of $(T^{(2)}U)^{*}$. Finally, note that $TM$ is a
subbundle of $T^{(2)}M$.

\begin{note}
{\rm If $f$ is 2-{\it flat at} $m$ then $t(f)(m) = 0$. We can also see
its dual approach. We then have a short exact sequence as given
below.
\begin{equation*}
0 \longrightarrow  F_{m}^{2}/F_{m}^{3} \xrightarrow{\ \ i \ \ }
F_{m}/F_{m}^{3} \xrightarrow{\ \ j \ \ }
F_{m}/F_{m}^{2} \longrightarrow 0.
\end{equation*}
$F_{m}/F_{m}^{2}$ is called the space of 1-jets at $m$ and is
denoted by $J_{m}$ (p.~376 of \cite{11}). This is
just the cotangent space $T^{*}_{m}M$ (also called 1-covelocity
space). Two functions $f, g \in F_{m}$ are equivalent in $J_{m}$
if in a chart $(U, \va)$ at $m$ (and hence in any chart),
$f_{\va}$ and $g_{\va}$ have the same first-order partial
derivatives at $\va m$.

$F_{m}/F_{m}^{3}$ is called the space of 2-jets at $m$ and is
denoted by $J_{m}^{2}$ (2-covelocity space). Two functions $f, g
\in F_{m}$ are equivalent in $J_{m}^{2}$ if the first and the
second partial dervatives of $f_{\va }$ and $g_{\va}$ agree at
$\va m$ in any chart $(U, \va)$ at $m$.

Let us now come to $F_{m}^{2}/F_{m}^{3}$. We show that this can be
identified with $S_{m}^{2}$, the space of symmetric bilinear
functionals on $T_{m}M\times T_{m}M$. Let ${\bf v, w} \in T_{m}M$.
Let $X$ and $Y$ be any two vector fields on $M$ such that $X(m) =
{\bf v}$ and $Y(m) = {\bf w}$. Let $f\in F_{m}^{2}$. Note that for
any vector field $Z$ on $M$, $Zf \in F_{m}$ and so $Zf(m) = 0$.
Now we have $X(Yf)(m) - Y(Xf)(m) = [X, Y]f(m) = 0$, so that
$X(Yf)(m) = Y(Xf)(m) = {\bf v}(Yf)(m) = {\bf w}(Xf)(m)$. This
shows that $D^{2}f(m)({\bf v, w}) = X(Yf)(m)$ is well-defined and
belongs to $S_{m}^{2}$. Also it may be checked that if $(U, \va )$
is any chart at $m$, then $D^{2}f(m)({\bf v, w}) =
D^{2}f_{\va}(\va m)({\bf v}_{\va}, {\bf w}_{\va})$,  where
$D^{2}f_{\va}$ denotes the second-order Frechet derivative of
$f_{\va }$ at $\va m$. $D^{2}f(m)$ is called the {\it Hessian} of
$f$ at $m$ at the critical point $m$. Clearly if $f \in F_{m}^{3}$
then $Df(m) = 0$ and $D^{2}f(m) = 0$. It follows that $F_{m}^{3}$
is the kernel of the mapping $f\mapsto D^{2}f(m)$.

Let us get back to our short exact sequence which we write as
\begin{equation*}
0 \longrightarrow  S_{m}^{2} \xrightarrow{\ \ i \ \ }
J_{m}^{2} \xrightarrow{\ \ j \ \ } J_{m} \longrightarrow
0.
\end{equation*}
Here $J_{m}$ is a vector space of dimension $d$, $J_{m}^{2}$ is a
vector space of dimension $d + \left(\begin{smallmatrix} d \\ 2\end{smallmatrix} \right)$ and $S_{m}^{2}$
is a vector space of dimension $\left( \begin{smallmatrix} d \\ 2\end{smallmatrix} \right)$. We can define
in the usual way the three vector bundles $J$, $J^{2}$ and $S^{2}$
on $M$ with fibers $J_{m}$, $J^{2}_{m}$ and $S^{2}_{m}$
respectively, to have a short exact sequence of vector bundles on
$M$:
\begin{equation*}
0 \longrightarrow  S^{2} \xrightarrow{\ \ i \ \ } J^{2}
\xrightarrow{\ \ j \ \ } J \longrightarrow 0.
\end{equation*}

It is natural to ask when this exact sequence splits. This leads
to the concept of a `dissection' as a splitting of this exact
sequence. Following this second-order Hessian structures and
dissections are in one-to-one correspondence \cite{3}.

Also we can define the `tangent bundle of order 2', denoted by $T^{2}M$
as a bundle of 2-jets \cite{11} defined in the following way, which is
certainly a different approach to $T^{(2)}M$ (the second-order tangent
bundle as defined in \S1.1).}
\end{note}

\subsection{\it The tangent bundle of order 2, \hbox{$T^{2}M$}}

Let $M$ be a $d$-dimensional manifold (p.~368 of \cite{11}). The
tangent bundle of order 2, $T^{2}M$, of $M$ is the $3d$-dimensional
manifold of 2-jets $j^{2}f$ at $0 \in {\bf R}$ of differentiable
curves $f\hbox{\rm :}\ {\bf R} \rightarrow M $. ($T^{2}M$ is also called
2-velocity space.)

$T^{2}M$ has a natural bundle structure over $M$, $\pi ^{2}\hbox{\rm :}\ s
T^{2}M \rightarrow M $ denotes the canonical projection. The
tangent bundle $TM$ is nothing but the manifold of 1-jets $j^{1}f$
at $0 \in {\bf R}$ of $f\hbox{\rm :}\ {\bf R} \rightarrow M $. If we denote
$\pi ^{12}\hbox{\rm :}\ T^{2}M \rightarrow TM $ the canonical projection,
then $T^{2}M$ has a bundle structure over $TM$. Note that $T^{2}M$
is {\it not} a vector bundle. (The bundle of 2-jets defined above.)
There is a result stating that `the linear connection $\nabla $
on $M$ determines a vector bundle structure on $\pi ^{2}\hbox{\rm :}\ T^{2}M
\rightarrow M $ and a vector bundle isomorphism $ T^{2}M
\rightarrow TM \bigoplus TM$'.

Note that the space $T^{2}M$ of 2-velocities on $M$ may be
identified with a submanifold of $T(TM)$, the tangent bundle to
$TM$ (see p.~372 of \cite{11}).\vspace{-.3pc}

\section{Hessian structures}

Let $\Omega ^{n}(M)$ denote the space of covariant tensors of
$n$th order on $M$. Since we wish to associate $n$th derivatives
to smooth real-valued functions on manifolds, we have the
following definition.

\setcounter{theore}{0}
\begin{definit}$\left.\right.$\vspace{.5pc}

\noindent {\rm An {\it $n$-th order Hessian structure} $H^{n}$ on $M$ is a mapping
$H^{n}\hbox{\rm :}\ {\cal F}(M)\rightarrow \Omega ^{n}(M)$, $f\mapsto
H^{n}f$, such that (i) $H^{n}$ is real linear in $f$ and (ii) if
$f \in F_{m}^{n}$ then $H^{n}f(m) = D^{n}f(m)$.}
\end{definit}

Note that the properties of $H^{n}$ imply that  if $f \in  {\cal
F}(M)$ is constant on  an  open  set $U$  in $M$,  then $H^{n}f$
is zero on $U$. Consequently, $H^{n}$ is localizable so that it is
meaningful to  talk of $H^{n}f$ for $f \in  {\cal F}(U)$. $H^{n}$
is said to be a symmetric Hessian structure if $H^{n}f$ is a
symmetric covariant tensor of order $n$, for all $f$. In \cite{3}, a
(second-order) Hessian structure was by definition symmetric. Here
we drop the requirement of symmetry from the definition.

Suppose now $\n $ is a connection on $M$. i.e., a connection on the
tangent bundle $TM$ of $M$. Recall \cite{13} that if $\omega \in \Omega
^{n}(M)$ then $\n \omega \in \Omega ^{n+1}(M)$ is defined by $ \n \omega
(X_{0}, \ldots, X_{n}) = (\n _{X_{0}} \omega )(X_{1}, \ldots, X_{n})$.
Here $(\n _{X_{0}}\omega )$ is defined by
\begin{align*}
(\n _{X_{0}} \omega )(X_{1}, \ldots, X_{n}) &= \n _{X_{0}} (\omega
(X_{1}, \ldots, X_{n}))\\
&\quad \, - \sum_{i=1}^{n}\omega (X_{1}, \ldots, \n _{X_{0}}X_{i}, \ldots, X_{n}).
\end{align*}
$\n \omega $ is called the {\it covariant derivative} of $\omega $ with
respect to the connection $\n $ on $M$. Also note the following:

\begin{enumerate}
\renewcommand\labelenumi{(\arabic{enumi})}
\item Let $\omega $
be a symmetric covariant tensor of order $n$. Then $\n \omega $ is
symmetric if and only if $(\n \omega )(X_{0}, X_{1}, \ldots, X_{n}) =
(\n \omega )(X_{1}, X_{0}, \ldots, X_{n})$, for all $X_{i} \in \sx(M)$.

\item Let $f \in {\cal F}(M)$. Define $\n f = df$, $\n ^{2}f = \n (df)$,
and recursively $\n ^{n+1} f = \n (\n ^{n} f)$. For $X, Y \in \sx(M)$,
$\n ^{2}f(X, Y) = X \. Y\. f -\n _{X}Y \. f$.
\end{enumerate}

\begin{theor}[\!]
Let $\n $ be a connection on $M$. Then $\n ^{n}f$ is an $n$-th
order Hessian structure on $M${\rm ,} for all $n \geq 1$.
\end{theor}

\begin{proof}
We first claim that if $f\in F_{m}^{n+k}$, then for fixed $X_{1}, \ldots,
X_{n}$, $\n ^{n}f(X_{1}, \ldots, X_{n}) = X_{1}\. \ldots \. X_{n}\.
f + g$, where $g\in F_{m}^{k+1}$. (This also implies that $\n
^{n}f(X_{1}, \ldots, X_{n})\in F_{m}^{k}$.) The proof is by induction
on $n$. Clearly it is true for $n=2$. Suppose it is true for $n$. Let
$f\in F_{m}^{n+1+k}$. Then
\begin{align*}
&\n ^{n+1}f(X_{0}, X_{1}, \ldots, X_{n}) \\
&\quad\, = [\n _{X_{0}}(\n ^{n}f)](X_{1}, \ldots, X_{n})\\
&\quad\, = \n _{X_{0}}[(\n ^{n}f(X_{1},\ldots, X_{n})]
 - \sum _{i=1}^{n}  \n ^{n}f(X_{1},\ldots, \nabla _{X_{0}}X_{i}, \ldots, X_{n})\\
&\quad\, = X_{0}(X_{1}\. X_{2}\. \cdots \. X_{n}\. f + g)
 - \sum _{i=1}^{n}\n ^{n}f(X_{1}, \ldots , \n_{X_{0}}X_{i}, \ldots, X_{n})\\
&\quad\, = X_{0}\. \cdots \. X_{n}\. f + X_{0}\. g
 - \sum _{i=1}^{n}\n ^{n}f(X_{1}, \ldots, \n_{X_{0}}X_{i}, \ldots, X_{n}).
\end{align*}
Here $g\in F_{m}^{k+2}$ and therefore $X_{0} g\in F_{m}^{k+1}$.
Also each term inside the summation on the right-hand side
belongs to ${\bf F}_{m}^{k+1}$. So the claim is proved.

Now we prove that $\n ^{n}f$ is an $n$th order Hessian structure
on $M$. Now the linearity in $f$ is trivially true and so we have
only to prove that if $f\in F_{m}^{n}$, then $\n ^{n}f(X_{1},
\ldots , X_{n})(m)=X_{1}\. \cdots \. X_{n} \cdot f(m)$. But our
claim implies that the difference $\n ^{n}f(X_{1}, \ldots , X_{n})
- X_{1}\. \cdots \. X_{n}\. f\in F_{m}$ if $f\in F_{m}^{n}$. \hfill q.e.d.
\end{proof}

In the case of ${\bf R}^{n}$, the $n$th derivative is symmetric.
In case of a manifold, we ask for conditions under which the $n$th
derivative is symmetric. We have the following results.

\begin{theor}[\!]$\left.\right.$\vspace{-.5pc}

\begin{enumerate}
\renewcommand\labelenumi{\rm (\alph{enumi})}
\leftskip .15pc
\item $\n ^{2}f$ is symmetric if and only if $\n$ has torsion zero.

\item $\n ^{2}f$ and $\n ^{3}f$ are both symmetric if and only if
$\n $ has both torsion and curvature\break zero.

\item If $\n $ has torsion and curvature zero then $\n ^{n}f$ is
symmetric{\rm ,} for all $n \geq 2$.\vspace{-.6pc}
\end{enumerate}
\end{theor}

\begin{proof}
Recall the definitions of torsion $T$ and curvature
$R$ from \cite{4} defined by $ T(X, Y) = \n _{X}Y - \n _{Y}X - [X, Y]$
and $R(X, Y)Z = \n _{X}\n _{Y}Z - \n _{Y}\n _{X}Z - \n _{[X,
Y]}Z.$ The proofs of (a) and (b) are trivial. (c) Assume that  $\n $ has torsion and curvature zero. We need to
prove that $\n ^{n}f$ is symmetric, for all $n \geq 2$. The proof
is by induction on $n$. By (a) and (b) the result is true for $n=
2, 3$. Suppose the result is true for $n$. We have to prove that
the result is true for $(n+1)$. Consider
\begin{align*}
&\n ^{n+1}f(X_{0}, X_{1}, \ldots, X_{n})\\[.4pc]
&\quad\, =  [\n _{X_{0}}(\n ^{n}f)](X_{1}, \ldots, X_{n})\\[.4pc]
&\quad\, = \n _{X_{0}}(\n ^{n}f(X_{1}, \ldots, X_{n})) -
     \sum_{i=1}^{n}\n ^{n}f (X_{1}, \ldots, \n _{X_{0}}X_{i}, \ldots, X_{n})\\[.4pc]
&\quad\, = \n _{X_{0}} \bigg[ \n _{X_{1}} (\n ^{n-1} f(X_{2}, \ldots, X_{n})) \\[.4pc]
&\qquad\, -  \sum_{j=2}^{n}\n ^{n-1}f (X_{2}, \ldots, \n _{X_{1}}X_{j}, \ldots, X_{n}) \bigg]\\[.4pc]
&\qquad \,-  \sum_{i=1}^{n}\n ^{n}f ( X_{1}, \ldots, \n _{X_{0}}X_{i}, \ldots, X_{n})\\[.4pc]
&\quad\, = \n _{X_{0}}  [\n _{X_{1}}(\n ^{n-1}f(X_{2}, \ldots,
X_{n})) ]\\[.4pc]
&\qquad\, - \sum_{j=2}^{n}\n _{X_{0}} (\n ^{n-1}f ( X_{2}, \ldots, \n _{X_{1}}X_{j}, \ldots, X_{n})\\[.4pc]
&\qquad\, - \bigg[ \n ^{n}f ( \n _{X_{0}}X_{1}, X_{2}, \ldots, X_{n})\\[.4pc]
&\qquad\, + \sum_{i=2}^{n}\n ^{n}f (X_{1}, X_{2}, \ldots, \n _{X_{0}}X_{i}, \ldots, X_{n})\bigg] \\[.4pc]
&\quad\, = \n _{X_{0}}\n _{X_{1}}(\n ^{n-1}f(X_{2}, \ldots, X_{n}))\\[.4pc]
&\qquad\, -  \sum_{j=2}^{n}\n _{X_{0}} (\n ^{n-1}f ( X_{2}, \ldots, \n _{X_{1}}X_{j}, \ldots, X_{n}))
\end{align*}
\begin{align*}
&\qquad\,- \bigg[ \n _{\n _{X_{0}}X_{1}}(\n ^{n-1}f(X_{2}, \ldots,
X_{n}))\\[.3pc]
&\qquad\, - \sum_{j=2}^{n}\n ^{n-1}f (X_{2}, \ldots, \n _{\n _{X_{0}}X_{1}}X_{j}, \ldots, X_{n}) \bigg]\\[.3pc]
&\qquad\, - \sum_{i=2}^{n}\n ^{n}f ( X_{1}, X_{2}, \ldots, \n _{X_{0}}X_{i}, \ldots, X_{n})\\[.3pc]
&\quad\, = [ \n _{X_{0}}\n _{X_{1}}(\n ^{n-1}f(X_{2}, \ldots, X_{n}))
- \n _{\n _{X_{0}}X_{1}}(\n ^{n-1}f(X_{2}, \ldots, X_{n})) ]\\[.3pc]
&\qquad\, - \sum_{j=2}^{n}\n _{X_{0}} (\n ^{n-1}f (X_{2}, \ldots,
 \n _{X_{1}}X_{j}, \ldots, X_{n}))\\[.3pc]
&\qquad\, + \sum_{j=2}^{n}\n ^{n-1}f (X_{2}, \ldots,\n _{\n
_{X_{0}}X_{1}}X_{j}, \ldots, X_{n})\\[.3pc]
&\qquad\,-  \sum_{i=2}^{n}\n ^{n}f (X_{1}, X_{2}, \ldots, \n
_{X_{0}}X_{i}, \ldots, X_{n}).
\end{align*}
Interchanging $X_{0}$ and $X_{1}$ in the above formula and using
$T(X, Y)\. f = \n _{X}Y\. f - \n _{Y}X\. f - [X, Y]\. f$, we
consider the difference
\begin{align}
&\n ^{n+1}f(X_{0}, X_{1}, \ldots, X_{n}) -
 \n ^{n+1}f(X_{1}, X_{0}, \ldots, X_{n})\nonumber\\[.3pc]
&\quad\, = T(X_{1}, X_{0})\. \n ^{n-1}f(X_{2}, \ldots, X_{n})\nonumber\\[.3pc]
&\qquad\, +  \l [ \sum_{j=2}^{n}\n _{X_{1}} (\n ^{n-1}f ( X_{2},
\ldots, \n _{X_{0}}X_{j}, \ldots, X_{n})) \nonumber \r.\\[.3pc]
&\qquad\, - \l. \sum_{j=2}^{n}\n _{X_{0}} (\n ^{n-1}f ( X_{2},
\ldots, \n _{X_{1}}X_{j}, \ldots, X_{n})) \r]\nonumber\\[.3pc]
&\qquad\, + \l [ \sum_{j=2}^{n}\n ^{n-1}f (X_{2}, \ldots, \n _{\n
_{X_{0}}X_{1} - \n _{X_{1}}X_{0}}X_{j}, \ldots, X_{n}) \r ]\nonumber\\[.3pc]
&\qquad\, + \bigg[ \sum_{i=2}^{n}\n ^{n}f (X_{0}, X_{2}, \ldots, \n _{X_{1}}X_{i},
 \ldots, X_{n})\nonumber\\[.3pc]
&\qquad\, -  \sum_{i=2}^{n}\n ^{n}f ( X_{1}, \ldots, \n _{X_{0}}X_{i},
 \ldots, X_{n}) \bigg].
\end{align}
Now we consider
\begin{align}
&\sum_{i=2}^{n}\n ^{n}f (X_{1}, \ldots, \n
_{X_{0}}X_{i}, \ldots, X_{n})\nonumber
\end{align}
\begin{align}
&\quad\, = \sum_{i=2}^{n} [\n _{X_{1}} (\n ^{n-1}f)]
 (X_{2}, \ldots, \n _{X_{0}}X_{i}, \ldots, X_{n})\nonumber\\[.3pc]
&\quad\, = \sum_{i=2}^{n} \bigg[\n _{X_{1}} (\n ^{n-1}f(X_{2}, \ldots, \n _{X_{0}}X_{i}, \ldots, X_{n})) \nonumber\\[.3pc]
&\qquad\, -  \sum_{k=2}^{n}\n ^{n-1}f (X_{2}, \ldots, \n
_{X_{0}}X_{i}, \ldots, \n _{X_{1}}\n_{X_{0}}X_{k}, \ldots, X_{n})
 \bigg].
\end{align}
Similarly we have,
\begin{align}
&\sum_{i=2}^{n}\n ^{n}f (X_{0},X_{2}, \ldots, \n
_{X_{1}}X_{i}, \ldots, X_{n})\nonumber\\[.25pc]
&\quad\, = \sum_{i=2}^{n} \bigg[
  \n _{X_{0}} (\n ^{n-1}f(X_{2}, \ldots, \n _{X_{1}}X_{i}, \ldots, X_{n})) \nonumber\\[.25pc]
&\qquad \, - \sum_{k=2}^{n}\n ^{n-1}f (X_{2}, \ldots, \n
_{X_{1}}X_{i}, \ldots, \n _{X_{0}}\n_{X_{1}}X_{k}, \ldots, X_{n}) \bigg].
\end{align}
Substitute (2) and (3) in (1) to get
\begin{align*}
&\n ^{n+1}f(X_{0}, X_{1}, \ldots, X_{n}) -
 \n ^{n+1}f(X_{1}, X_{0}, X_{2}, \ldots, X_{n})\\[.25pc]
&\quad\,=  T(X_{1}, X_{0})\. \n ^{n-1}f(X_{2}, \ldots, X_{n})\\[.25pc]
&\qquad\, + \sum_{j=2}^{n}\n ^{n-1}f (X_{2}, \ldots,
 \n _{\n _{X_{0}}X_{1} - \n _{X_{1}}X_{0}}X_{j}, \ldots, X_{n})\\[.25pc]
&\qquad\, + \sum_{i=2}^{n}\sum_{k=2}^{n} \n ^{n-1}f ( X_{2}, \ldots,
\n _{X_{0}}X_{i}, \ldots, \n _{X_{1}}\n _{X_{0}}X_{k}, \ldots, X_{n})\\[.25pc]
&\qquad\, - \sum_{i=2}^{n}\sum_{k=2}^{n} \n ^{n-1}f (X_{2}, \ldots,
\n _{X_{1}}X_{i}, \ldots, \n _{X_{0}}\n _{X_{1}}X_{k}, \ldots, X_{n})\\[.25pc]
&\quad\,= T(X_{1}, X_{0})\. \n ^{n-1}f(X_{2}, \ldots, X_{n})\\[.25pc]
&\qquad\, + \sum_{j=2}^{n}\n ^{n-1}f ( X_{2}, \ldots,
 \n _{T(X_{0}, X_{1})}X_{j}, \ldots, X_{n})\\[.25pc]
&\qquad\, - \sum_{i=2}^{n}\n ^{n-1}f (X_{2}, \ldots, R(X_{0},
X_{1})X_{i}, \ldots, X_{n})
\end{align*}
This is enough to prove the theorem. \hfill q.e.d.
\end{proof}

Ambrose, Palais and Singer in \cite{3} proved that every symmetric
second-order Hessian structure $H^{2}$ is of the form $\n ^{2}$,
for a symmetric connection $\n $. We see below that the condition
of symmetry is unnecessary. Before this, we prove the following Lemma~2.4.
\newpage

\begin{lem}
Fix ${\bf v, w} \in T_{m}M$. Let $Bf(m)({\bf v,w})=D^{2}f(m)({\bf
v, w}) - Hf(m)({\bf v,w})$. Then
\begin{enumerate}
\renewcommand\labelenumi{\rm (\roman{enumi})}
\leftskip .35pc
\item the map $B\hbox{\rm :}\ f \rightarrow Bf(m)({\bf v,w})$ is
real linear in $f${\rm ,}

\item if $df (m) = 0$ then $Bf(m)=0${\rm ,}

\item $B$ is a derivation in $f${\rm ,} i.e.{\rm ,} $B(fg)=fB(g) + gB(f)$  for
all $f,g \in {\cal F}(M)$.
\end{enumerate}
\end{lem}

\begin{proof}
Lemma 2.4 (i) and (ii) follow from the relevant definitions.
Lemma~2.4 (iii) follows from Lemma~2.4 (ii) by noting that if we
set $h=(f-fm)(g-gm)$, then $dh(m)=0$ so that $Bh(m)=0$.\hfill {\rm
q.e.d.}
\end{proof}

\begin{theor}[\!]
If $H^{2}$ is a second-order Hessian structure on $M${\rm ,} then
there exists a unique connection $\n $ on $M$ such that $H^{2}f =
\n ^{2}f$ for all $f$. Thus there is a one-to-one correspondence
between second-order Hessian structures on $M$ and connections on
$M$.
\end{theor}

\begin{proof}
Suppose $H^{2}$ is a second-order Hessian structure on $M$. Fix
$X, Y \in \sx(M)$ and $f \in {\cal F}(M)$. Define $Bf(X, Y) = X\.
Y \. f - H^{2}f(X, Y)$. Then $B\hbox{\rm :}\ f \mapsto Bf(X, Y)$
is a derivation. Hence it defines a vector field on $M$ which we
denote it by $\n _{X}Y$ such that $Bf(X, Y) = \n _{X}Y \. f $.
Then we have $\n _{X}Y \. f = X\. Y\. f - H^{2}f(X, Y)$. It can be
easily checked that $\n $ satisfies all the properties of a
connection and that  $\n $ is symmetric if $H^{2}$ is symmetric.
\hfill {\rm q.e.d.}
\end{proof}

It is now natural to ask the following question. Does every
third-order Hessian structure $H^{3}$ arise as $\n ^{3}$, for a
connection $\n $? We show below that this need not be the case.
However, every $H^{3}$ comes from a connection on the second-order
tangent bundle $T^{(2)}M$ of $M$. We characterize the connections
on $T^{(2)}M$ which arise in this way. We also characterize the
smaller class of connections which are associated to $H^{3}$ of
the form $\n ^{3}$. Incidentally, we prove that every connection
$\n $ on $TM$ induces a connection $\widetilde{\n}$ on $T^{(2)}M$.
We do all this in the language of charts and Frechet calculus
methods.

\section{Second-order connections}

We define vector bundles now, to recall notation \cite{2}. Let
{\bf E} and {\bf F} be finite dimensional Banach spaces with $U$
open in {\bf E}. We call the product $U \times {\bf F}$ a local
vector bundle (l.v.b.). We call $U$ the base space, which can be
identified with $U\times \{0\}$, which is called  the zero
section. For $u \in U$, $\{u\}\times {\bf F}$ is called the fiber
over $u$, which can be endowed with the Banach space structure of
{\bf F}. The map $\pi\hbox{\rm :}\ U\times {\bf F}\rightarrow U$
given by $\pi (u, x) = u$ is called the projection of $U\times
{\bf F}$. Note that $U \times {\bf F}$ is an open subset of ${\bf
E}\times {\bf F}$ and so is a local manifold. Let $U\times {\bf
F}$ and $U'\times {\bf F}'$ be two l.v.b.'s. A map
$\varphi\hbox{\rm :}\ U\times {\bf F}\rightarrow U'\times {\bf
F}'$ is called a l.v.b. mapping if $\varphi$ is smooth and it has
the form $\varphi (u, x) = (\varphi_{1}(u), \varphi_{2}(u)\cdot
x)$ where $\varphi_{1}\hbox{\rm :}\ U\rightarrow U'$ and
$\varphi_{2}\hbox{\rm :}\ U\rightarrow L({\bf F, F'})$ are smooth.
A l.v.b. map that is a diffeomorphism is called a l.v.b.
isomorphism. For example, (i) any linear map $A\in L({\bf E, F})$
defines a vector bundle map $\varphi_{A}\hbox{\rm :}\ {\bf
E}\times {\bf E}\rightarrow {\bf E}\times {\bf F}$ by
$\varphi_{A}(u, {\bf e})= (u, A\cdot {\bf e})$ and (ii) if $U$ is
open in {\bf E}, $V$ is open in {\bf F} and $f\hbox{\rm :}\
U\rightarrow V$ is smooth, then the  map $Tf\hbox{\rm :}\ U\times
{\bf E}\rightarrow V\times {\bf F}$ given by $Tf(u, {\bf e})=
(f(u), Df(u)\cdot {\bf e})$ is a l.v.b. map.

In this setting, let $\pi\hbox{\rm :}\ E \rightarrow M $ be a
vector bundle, where $M$ is a manifold modelled on a finite
dimensional Banach space {\bf E} and each fibre of $E$ is modelled
on a finite dimensional Banach space {\bf F}. For the local
structure of vector bundles, we follow the notation given below.
Let $m\in M$. If $(U, \varphi )$ is a local chart at $m$, then the
vector bundle chart is a triple $(\Phi, \varphi, U)$ where the
following diagram commutes:
\newpage

$\left.\right.$\vspace{-1.75cm}

\begin{center}
\unitlength=.09mm
\begin{picture}(500,500)(0,0)
\put(140,140){$U$} \put(300,150){$\varphi U$}
\put(40,295){$\pi^{-1}(U)$} \put(300,300){$\varphi U\times {\bf
F}$} \put(150,285){\vector(0,-1){110}}
\put(300,285){\vector(0,-1){110}} \put(175,300){\vector(1,0){110}}
\put(175,150){\vector(1,0){110}} \put(225,125){$\varphi$}
\put(225,310){$\Phi$} \put(115,225){$\pi$}
\put(310,225){$Pr_{1}$\qquad.}
\end{picture}
\vspace{-.9cm}
\end{center}

Thus, if $\tau\hbox{\rm :}\ U\rightarrow E$ is a local section of
$\pi$, then the principal part $\tau_{\varphi}\hbox{\rm :}\
\varphi U\rightarrow {\bf F}$ with respect to the vector bundle
chart $(\Phi, \varphi, U)$ is given by $(\Phi \circ \tau \circ
\varphi^{-1})(u)=(u, \tau_{\varphi}(u))$; ($u\in \varphi U$).
Similarly, for a vector field $X$ of $M$, the local representative
$X_{\varphi}\hbox{\rm :}\ \varphi U\rightarrow {\bf E}$ with
respect to the chart $(U, \varphi)$ is given by $(T\varphi \circ X
\circ \varphi^{-1})(u)= (u, X_{\varphi}(u))$; ($u\in \varphi U$).
For any two vector bundle charts $(\Phi, \varphi, U)$ and $(\Psi,
\psi, V)$ at $e \in E$, we have the l.v.b. mapping given by $(\Psi
\circ \Phi^{-1})(u, x) = (\varphi_{1}(u), \varphi_{2}(u)\cdot x)$
for all $u\in \varphi U$, $x\in {\bf F}$ where $\varphi_{1} =\psi
\circ \varphi^{-1}$ and $\varphi_{2}\hbox{\rm :}\ \varphi
U\rightarrow L({\bf F, F})$ given by $\va _{2}(u) \. x = Pr_{2}
\circ \Psi \circ \Phi ^{-1} \circ (u, x) \in {\bf F}$ is smooth.
Let $\sx_{E}(M)$ denote the smooth sections of the vector bundle
$\pi\hbox{\rm :}\ E\rightarrow M$ and let $\sx (M)$ denote the set
of all smooth vector fields on $M$. Following \cite{4} and
\cite{8} we have the definition given below.

\setcounter{theore}{0}
\begin{definit}$\left.\right.$\vspace{.5pc}

\noindent {\rm A {\it connection} $\widetilde{\nabla}$ on a vector bundle
$\pi\hbox{\rm :}\ E\rightarrow M$ is a mapping
$\widetilde{\nabla}\hbox{\rm :}\ \sx(M)\times
\sx_{E}(M)\rightarrow \sx_{E}(M)${\rm ,} $(X, \tau )\mapsto
\widetilde{\nabla}_{X}\tau $ for all $X\in \sx(M)${\rm ,} $\tau
\in \sx_{E}(M)$ such that it satisfies the following properties:
\begin{enumerate}
\renewcommand\labelenumi{\rm (\roman{enumi})}
\leftskip .1pc
\item $\widetilde{\n }$ is real bilinear in $X$ and $\tau${\rm ,}

\item $\widetilde{\n }_{f X}\tau = f\widetilde{\n }_{X} \tau$ and
$\widetilde{\n }_{X}(f \tau ) = X(f) \tau + f
\widetilde{\n}_{X}\tau ${\rm ,} for $f \in {\cal F}(M)$ and $\tau
\in \sx_{E}(M)$.
\end{enumerate}}
\end{definit}

Note that in a local vector bundle chart $(\Phi,\varphi,U)$ of
$E$, a connection  has the form $ (\widetilde{\nabla}_{X}\tau
)_{\varphi m} = D\tau_{\varphi m}\cdot X_{\varphi m} -
\widetilde{\Gamma}_{\varphi}(\varphi m)(X_{\varphi m},
\tau_{\varphi m})$ for all $m\in U$  where
$\widetilde{\Gamma}_{\varphi}\hbox{\rm :}\ \varphi U \rightarrow
L^{2}({\bf E}\times {\bf F}, {\bf F})$ is the Christoffel symbol
of $\widetilde{\nabla}$.

Let us consider the vector bundle $E=T^{(2)}M$ over $M$ (the
second-order tangent bundle of $M$). Locally, if $M=U'$ is open in
{\bf E}, then $T^{(2)}U'= U'\times {\bf E} \oplus ({\bf E}\Delta
{\bf E})$. We define a connection $\widetilde{\nabla }$ on
$T^{(2)}M$ as above and we now find the transformation property of
$\widetilde{\Gamma} $. Let $(U, \varphi)$ and $(V, \psi)$ be two
charts at $m$ and suppose $(\Phi, \varphi, U)$ and $(\Psi, \psi,
V)$ are vector bundle charts at $e\in E$ and let $F=\psi \circ
\varphi^{-1}$, $\varphi m= u$, $\psi m =   v$. Then $F(u)=v$. If
$\tau \hbox{\rm :}\  U\rightarrow T^{(2)}M$ is a local section of
$T^{(2)}M$, then the principal part is $\tau_{\varphi}\hbox{\rm
:}\ \varphi U \rightarrow {\bf F}$ with respect to the vector
bundle chart $(\Phi, \varphi, U)$, where ${\bf F} = {\bf E}\oplus
({\bf E}\Delta {\bf E})$. Then we have the following rules. We
write $\tau_{\varphi} = V_{\varphi}\oplus S_{\varphi}$, where
$V_{\varphi}(u) \in {\bf E}$ and $S_{\varphi}(u) \in {\bf E}\Delta
{\bf E}$. $ \tau_{\psi}(v) = \varphi_{2}(u)\cdot
\tau_{\varphi}(u)$ and $X_{\psi}(v) = DF(u) \cdot X_{\varphi}(u)$.

Then by the definition of second-order tangent vectors, we have
\begin{equation}
V_{\psi}(v) = DF(u)\cdot V_{\varphi}(u) + D^{2}F(u) \cdot
S_{\varphi}(u)
\end{equation}
and
\begin{equation}
S_{\psi}(v) = [DF(u)\Delta DF(u)]\cdot S_{\varphi}(u).
\end{equation}
From the definition of $\widetilde{\nabla}$, we have in a local
chart $(U, \varphi)$ of $M$,
\begin{equation*}
(\widetilde{\nabla}_{X}\tau )_{\varphi} = D\tau_{\varphi} \cdot
X_{\varphi}
  - \widetilde{\Gamma}_{\varphi}(X_{\varphi}, \tau_{\varphi}).
\end{equation*}
We write $ \widetilde{\Gamma}_{\varphi}(X_{\varphi},
\tau_{\varphi}) = \widetilde{\Gamma}^{1}_{\varphi}(X_{\varphi},
\tau_{\varphi}) \oplus
\widetilde{\Gamma}_{\varphi}^{2}(X_{\varphi}, \tau_{\varphi}), $
where  $\widetilde{\Gamma}_{\varphi}^{1}(X_{\varphi},
\tau_{\varphi})(u) \in {\bf E}$ and
$\widetilde{\Gamma}_{\varphi}^{2}(X_{\varphi}, \tau_{\varphi})(u)
\in {\bf E}\Delta {\bf E}$. Therefore
\begin{equation*}
(\widetilde{\nabla}_{X}\tau )_{\varphi} = [DV_{\varphi} \cdot
X_{\varphi}
  - \widetilde{\Gamma}_{\varphi}^{1}(X_{\varphi}, \tau_{\varphi})]\oplus
[DS_{\varphi} \cdot X_{\varphi}
  - \widetilde{\Gamma}_{\varphi}^{2}(X_{\varphi}, \tau_{\varphi})] \in {\bf F}.
\end{equation*}
But
\begin{align*}
(\widetilde{\nabla}_{X}\tau )_{\psi}
&= \varphi_{2}(\widetilde{\nabla}_{X}\tau )_{\varphi} \\[.3pc]
&= DF\cdot [DV_{\varphi} \cdot X_{\varphi}
       -\widetilde{\Gamma}_{\varphi}^{1}(X_{\varphi},
       \tau_{\varphi})]\\[.3pc]
       &\quad\, +
         D^{2}F\cdot [DS_{\varphi} \cdot X_{\varphi}
         - \widetilde{\Gamma}_{\varphi}^{2}(X_{\varphi},
         \tau_{\varphi})]\\[.3pc]
&\quad\,\oplus (DF\Delta DF)[DS_{\varphi} \cdot X_{\varphi}
  - \widetilde{\Gamma}_{\varphi}^{2}(X_{\varphi}, \tau_{\varphi})].
\end{align*}
On the other hand, $ (\widetilde{\nabla}_{X}\tau )_{\psi} =
D\tau_{\psi} \cdot X_{\psi} - \widetilde{\Gamma}_{\psi}(X_{\psi},
\tau_{\psi})$.  Now
\begin{align*}
[DV_{\psi}\cdot X_{\psi}](v)
&= DV_{\psi}(v) \cdot DF(u) \cdot X_{\varphi}(u)\\[.3pc]
&= D(V_{\psi} \circ F)(u) \cdot X_{\varphi }(u)\\[.3pc]
&= D[DF \cdot V_{\varphi} + D^{2}F \cdot S_{\varphi}](u) \cdot
X_{\varphi}(u) \ \ (\hbox{by eqs~(4) and
(5)})\\[.3pc]
&= D^{2}F(u)(X_{\varphi}(u)\Delta V_{\varphi}(u))
+ DF(u) \cdot DV_{\varphi}(u)\cdot X_{\varphi}(u)\\[.3pc]
&\quad\,+ D^{3}F(u)(X_{\varphi}(u)\Delta S_{\varphi}(u)) +
D^{2}F(u) \cdot DS_{\varphi}(u) \cdot X_{\varphi}(u).
\end{align*}
Again,
\begin{align*}
[DS_{\psi}\cdot X_{\psi}](v)
&= DS_{\psi}(v)  \cdot DF(u) \cdot X_{\varphi}(u)\\[.3pc]
&= D[S_{\psi}\circ F](u) \cdot X_{\varphi}(u)\\[.3pc]
&= D[(DF\Delta DF) \cdot S_{\varphi}](u) \cdot X_{\varphi}(u)\\[.3pc]
&= D(DF\Delta DF)(u)(X_{\varphi}(u)\Delta S_{\varphi}(u))\\[.3pc]
&\quad\, + (DF\Delta DF)(u) \cdot DS_{\varphi}(u) \cdot
X_{\varphi}(u).
\end{align*}
Therefore
\begin{align*}
&[DV_{\psi}\cdot X_{\psi}](v)\oplus [DS_{\psi}\cdot X_{\psi}](v)\\[.3pc]
&\quad\, = [D^{2}F(u)(X_{\varphi}(u)\Delta V_{\varphi}(u)) + DF(u) \cdot DV_{\varphi}(u)\cdot X_{\varphi}(u)\\[.3pc]
&\qquad\, + D^{3}F(u)(X_{\varphi}(u)\Delta S_{\varphi}(u)) + D^{2}F(u) \cdot DS_{\varphi}(u) \cdot X_{\varphi}(u)]\\[.3pc]
&\qquad\, \oplus \{ D(DF\Delta DF)(u)(X_{\varphi}(u)\Delta
S_{\varphi}(u))\\[.3pc]
&\qquad\, + (DF\Delta DF)(u) \cdot DS_{\varphi}(u) \cdot X_{\varphi}(u)\}\\[.3pc]
&\quad\, = \{[DF(u) \cdot DV_{\varphi}(u)\cdot X_{\varphi}(u)
 + D^{2}F(u) \cdot DS_{\varphi}(u) \cdot X_{\varphi}(u)]\\[.3pc]
&\qquad\,\oplus [(DF\Delta DF)(u) \cdot DS_{\varphi}(u) \cdot
X_{\varphi}(u) ]\} + \{[ D^{3}F(u)(X_{\varphi}(u)\Delta S_{\varphi}(u))\\[.3pc]
&\qquad\,+ D^{2}F(u)(X_{\varphi}(u)\Delta V_{\varphi}(u))] \oplus
[D(DF\Delta DF)(u)(X_{\varphi}(u)\Delta S_{\varphi}(u))]\}
\end{align*}
\begin{align*}
&\quad\,= [\varphi_{2}(u)\cdot  [DV_{\varphi}(u)\cdot
X_{\varphi}(u)\oplus DS_{\varphi}(u)\cdot X_{\varphi}(u) ] ]\\[.3pc]
&\qquad\,+ \{ [ D^{3}F(u)(X_{\varphi}(u)\Delta S_{\varphi}(u)) +
D^{2}F(u)(X_{\varphi}(u)\Delta V_{\varphi}(u))]\\[.3pc]
&\qquad\,\oplus [ D(DF\Delta DF)(u)(X_{\varphi}(u)\Delta
S_{\varphi}(u)) ] \}.
\end{align*}
Again\vspace{-.7pc}
\begin{align*}
&\widetilde{\Gamma}_{\psi}(X_{\psi}, \tau_{\psi})\\[.3pc]
&\quad\, = D\tau_{\psi}\cdot X_{\psi} -(\widetilde{\nabla}_{X}\tau )_{\psi}\\[.3pc]
&\quad\,= [DV_{\psi}\cdot X_{\psi}\oplus DS_{\psi}\cdot X_{\psi}] -(\widetilde{\nabla}_{X}\tau )_{\psi}\\[.3pc]
&\quad\,= [\varphi_{2}[DV_{\psi}\cdot X_{\psi}\oplus
DS_{\psi}\cdot X_{\psi}]] +[D^{3}F(X_{\varphi}\Delta S_{\varphi}) + D^{2}F(X_{\varphi}\Delta V_{\varphi})\\[.3pc]
&\qquad\,\oplus  D(DF\Delta DF)(X_{\varphi} \Delta S_{\varphi})]
-[\varphi_{2}\cdot [D\tau_{\varphi}\cdot X_{\varphi}
- \widetilde{\Gamma}_{\varphi}(\varphi)(X_{\varphi}, \tau_{\varphi})]]\\[.3pc]
&\quad\,= \varphi_{2}\cdot [\widetilde{\Gamma}_{\varphi}
(X_{\varphi}, \tau_{\varphi})] + [D^{3}F(X_{\varphi}\Delta
S_{\varphi}) + D^{2}F(X_{\varphi}\Delta V_{\varphi})\\[.3pc]
&\qquad\, \oplus  D(DF\Delta DF)(X_{\varphi} \Delta S_{\varphi})]\\[.3pc]
&\quad\,= [DF\cdot \widetilde{\Gamma}_{\varphi}^{1}(X_{\varphi},
\tau_{\varphi}) + D^{2}F\cdot \widetilde{\Gamma}_{\varphi}^{2}
(X_{\varphi}, \tau_{\varphi}) \oplus (DF \Delta DF)\cdot
\widetilde{\Gamma}_{\varphi}^{2} (X_{\varphi}, \tau_{\varphi})]\\[.3pc]
&\qquad\,+ [D^{3}F(X_{\varphi}\Delta S_{\varphi})
+D^{2}F(X_{\varphi}\Delta V_{\varphi})\oplus D(DF\Delta
DF)(X_{\varphi} \Delta S_{\varphi})]\\[.3pc]
&\quad\,= [DF\cdot \widetilde{\Gamma}_{\varphi}^{1}(X_{\varphi},
\tau_{\varphi}) + D^{2}F\cdot \widetilde{\Gamma}_{\varphi}^{2}
(X_{\varphi}, \tau_{\varphi}) + D^{3}F(X_{\varphi}\Delta S_{\varphi})\\[.3pc]
&\qquad\,+ D^{2}F(X_{\varphi}\Delta V_{\varphi})] \oplus [(DF
\Delta DF)\cdot \widetilde{\Gamma}_{\varphi}^{2}(X_{\varphi},
\tau_{\varphi})\\[.3pc]
&\qquad\, + D(DF\Delta DF)(X_{\varphi} \Delta S_{\varphi})].
\end{align*}

This is the transformation formula for $\widetilde{\Gamma}$. In
what follows we simplify the presentation of this formula. As
usual we write $\tau_{\varphi} = V_{\varphi}\oplus S_{\varphi}$,
where $V_{\varphi}(u) \in {\bf E}$ and $S_{\varphi}(u) \in {\bf
E}\Delta {\bf E}$, $u \in \va U$. We define\vspace{-.5pc}
\begin{align*}
A_{\va}(X_{\va}, V_{\va}) &= Pr_{1}\circ \widetilde{\G}_{\va}(X_{\va}, V_{\va}\oplus 0),\\[.3pc]
R_{\va}(X_{\va}, S_{\va}) &= Pr_{1}\circ \widetilde{\G}_{\va}(X_{\va}, 0 \oplus S_{\va}),\\[.3pc]
B_{\va}(X_{\va}, V_{\va}) &= Pr_{2}\circ \widetilde{\G}_{\va}(X_{\va}, V_{\va}\oplus 0),\\[.3pc]
C_{\va}(X_{\va}, S_{\va}) &= Pr_{2}\circ \widetilde{\G}_{\va}(X_{\va}, 0 \oplus S_{\va}),
\end{align*}
where $A_{\va}(X_{\va}, . )\hbox{\rm :}\  {\bf E} \rightarrow
{\bf E}$, $R_{\va}(X_{\va}, . ) \hbox{\rm :}\  {\bf E}\Delta {\bf
E} \rightarrow {\bf E}$, $B_{\va}(X_{\va }, . ) \hbox{\rm :}\ {\bf
E} \rightarrow {\bf E}\Delta {\bf E}$ and $C_{\va}(X_{\va}, . )
\hbox{\rm :}\ {\bf E}\Delta {\bf E} \rightarrow {\bf E}\Delta {\bf
E}$. Now we can represent $\widetilde{\G }$ in the matrix form as\vspace{-.3pc}
\begin{equation*}
\widetilde{\Gamma}_{\varphi}(X_{\varphi}, \tau_{\varphi}) =
\begin{pmatrix} A_{\varphi}(X_{\varphi}, V_{\varphi}) &R_{\varphi}(X_{\varphi}, S_{\varphi})\\[.3pc]
B_{\varphi}(X_{\varphi}, V_{\varphi}) &C_{\varphi}(X_{\varphi}, S_{\varphi})\end{pmatrix},
\end{equation*}
i.e., we can write $ \widetilde{\Gamma}_{\varphi}^{1}(X_{\varphi},
\tau_{\varphi}) = A_{\varphi}(X_{\varphi},
V_{\varphi})+R_{\varphi}(X_{\varphi}, S_{\varphi}) $ and $
\widetilde{\Gamma}_{\varphi}^{2}(X_{\varphi}, \tau_{\varphi}) =
B_{\varphi}(X_{\varphi}, V_{\varphi})+C_{\varphi}(X_{\varphi},
S_{\varphi}). $ Therefore
\pagebreak
\begin{align*}
&[A_{\psi }(X_{\psi }, V_{\psi }) + R_{\psi}(X_{\psi},
S_{\psi})]\oplus [B_{\psi}(X_{\psi}, V_{\psi})+C_{\psi}(X_{\psi}, S_{\psi})]\\[.3pc]
&\quad\, = [DF\cdot (A_{\varphi}(X_{\varphi}, V_{\varphi})
+ R_{\varphi}(X_{\varphi}, S_{\varphi}))\\[.3pc]
&\qquad\, + D^{2}F\cdot (B_{\varphi}(X_{\varphi}, V_{\varphi})
+ C_{\varphi}(X_{\varphi}, S_{\varphi})) \\[.3pc]
&\qquad\, + D^{3}F(X_{\varphi}\Delta S_{\varphi})
+ D^{2}F (X_{\varphi}\Delta V_{\varphi})]\\[.3pc]
&\qquad\, \oplus [(DF\Delta DF)\cdot (B_{\varphi}(X_{\varphi},
V_{\varphi}) + C_{\va}(X_{\va}, S_{\va}))\\[.3pc]
&\qquad\, + D(DF\Delta DF)(X_{\varphi}\Delta S_{\varphi})].
\end{align*}
Hence
\begin{align*}
A_{\psi}(X_{\psi}, V_{\psi})+R_{\psi}(X_{\psi}, S_{\psi})
&= DF\cdot [A_{\varphi}(X_{\varphi}, V_{\varphi})
    + R_{\varphi}(X_{\varphi}, S_{\varphi})]\\[.3pc]
    &\quad\,+ D^{2}F\cdot [B_{\varphi}(X_{\varphi}, V_{\varphi})
    + C_{\varphi}(X_{\varphi}, S_{\varphi})]\\[.3pc]
    &\quad\,+ D^{3}F(X_{\varphi}\Delta S_{\varphi})
    + D^{2}F(X_{\varphi}\Delta V_{\varphi})
\end{align*}
and
\begin{align*}
B_{\psi}(X_{\psi}, V_{\psi}) + C_{\psi}(X_{\psi}, S_{\psi})
&= (DF\Delta DF)\cdot (B_{\varphi}(X_{\varphi}, V_{\varphi})\\[.3pc]
&\quad\, + C_{\varphi}(X_{\varphi}, S_{\varphi})) + D(DF\Delta
DF)(X_{\varphi}\Delta S_{\varphi}).
\end{align*}
We now have the following formulas for $A_{\varphi}$, $B_{\varphi}$,
$R_{\varphi}$ and $C_{\varphi}$.

\begin{case}
{\rm If $S_{\varphi}=0$, then $S_{\psi}=0$ and
$V_{\psi}=DF\cdot V_{\varphi}$. Therefore from  above, we have
\begin{equation}
A_{\psi}(X_{\psi}, V_{\psi}) = DF\cdot A_{\varphi}(X_{\varphi},
V_{\varphi}) + D^{2}F\cdot B_{\varphi}(X_{\varphi}, V_{\varphi}) + D^{2}F(X_{\varphi}\Delta V_{\varphi})
\end{equation}
and
\begin{equation}
B_{\psi}(X_{\psi}, V_{\psi}) = (DF\Delta DF)\cdot B_{\varphi}(X_{\varphi}, V_{\varphi}).
\end{equation}}
\end{case}

\begin{case}
{\rm If $V_{\varphi}=0$ then $V_{\psi}=D^{2}F\cdot
S_{\varphi}$. Then from  above, we have
\begin{align*}
R_{\psi}(X_{\psi}, S_{\psi}) &= DF\cdot R_{\varphi}(X_{\varphi}, S_{\varphi})\\[.3pc]
     &\quad\, + D^{2}F\cdot C_{\varphi}(X_{\varphi}, S_{\varphi}) + D^{3}F(X_{\varphi}\Delta S_{\varphi})\\[.3pc]
     &\quad\, - A_{\psi}(X_{\psi}, D^{2}F \cdot S_{\varphi})
\end{align*}
and
\begin{align*}
C_{\psi}(X_{\psi}, S_{\psi}) &= (DF\Delta DF)\cdot C_{\varphi}(X_{\varphi}, S_{\varphi})\\
&\quad\, + D(DF\Delta DF)(X_{\varphi}\Delta S_{\varphi}) - B_{\psi}(X_{\psi}, D^{2}F \cdot S_{\varphi}).
\end{align*}}
\end{case}

\begin{note}
{\rm If we choose $A_{\va}(X_{\va}, V_{\va}) = 0$ and
$B_{\va}(X_{\va}, V_{\va}) = - X_{\va} \Delta V_{\va}$,
where $\tau _{\va} = V_{\va} \oplus S_{\va}$, $V_{\va}(u)\in {\bf
E}$, $S_{\va}(u) \in {\bf E}\Delta {\bf E}$ ($u\in U$), then
$A_{\psi}(X_{\psi}, V_{\psi}) = 0$ and
$B_{\psi}(X_{\psi}, V_{\psi}) = - X_{\psi} \Delta V_{\psi}$.}
\end{note}

Using equations (6) and (7), we can easily see that the choice of
$A_{\va}$ and $B_{\va}$ are independent of charts and hence they
define a subclass of second-order connections on $T^{(2)}M$. We
have the following definition.

\begin{definit}$\left.\right.$\vspace{.5pc}

\noindent {\rm A connection $\widetilde{\nabla}$ on the vector
bundle $T^{(2)}M$ is called a {\it special second-order
connection} on $M$ if the Christoffel symbols of
$\widetilde{\nabla}$ in the given chart $(U, \varphi)$ are given
by the matrix of symbols
\begin{equation*}
\widetilde{\Gamma}_{\varphi}(X_{\varphi}, \tau_{\varphi})=
\begin{pmatrix}
0 &R_{\varphi}(X_{\varphi}, S_{\varphi})\\[.5pc]
-X_{\varphi}\Delta V_{\varphi} &C_{\varphi}(X_{\varphi}, S_{\varphi})
\end{pmatrix},
\end{equation*}
where $\tau_{\varphi} = V_{\varphi} \oplus S_{\varphi}$,
$V_{\varphi}(u) \in {\bf E}$, $S_{\varphi}(u) \in {\bf E}\Delta
{\bf E}$ ($u\in \varphi U$). We call $R_{\va}$ and $C_{\va}$ as
the Christoffel symbols of $\widetilde{\n}$ in the chart $(U,
\va)$.}
\end{definit}

\begin{note}
{\rm Suppose $\widetilde{\nabla}$ is a special second-order
connection on $M$. Then in any other chart $(V, \psi)$ at $m \in
M$, the transformation laws for the Christoffel symbols are given
by
\begin{equation}
R_{\psi}(X_{\psi}, S_{\psi}) = DF\cdot R_{\varphi}(X_{\varphi},
S_{\varphi}) + D^{2}F\cdot C_{\varphi}(X_{\varphi}, S_{\varphi}) +
D^{3}F(X_{\varphi}\Delta S_{\varphi})
\end{equation}
and
\begin{align}
C_{\psi}(X_{\psi}, S_{\psi}) &= (DF\Delta DF)\cdot
C_{\varphi}(X_{\varphi}\Delta S_{\varphi})
  + D(DF\Delta DF)(X_{\varphi}\Delta S_{\varphi})\nonumber\\[.3pc]
  &\quad\, + X_{\psi} \Delta  D^{2}F \cdot S_{\varphi}.
\end{align}
If $S_{\varphi } = Y_{\varphi }\Delta Z_{\varphi }$, we write
$R_{\varphi }(X_{\varphi} ,Y_{\varphi }, Z_{\varphi }) =R_{\varphi
}(X_{\varphi }, Y_{\varphi }\Delta Z_{\varphi })$ and similarly
for $C_{\varphi }$.}\vspace{-.4pc}
\end{note}

\section{Third-order Hessian structures}

Let $K$ be a third-order Hessian structure on $M$ (the case $n=3$
of Definition 2.1). In what follows we state that the following
lemma and the proof is easy and hence left to the\break reader.

\setcounter{theore}{0}
\begin{lem}
Let $F \hbox{\rm :}\  U \rightarrow V$ be a smooth map{\rm ,} where $U,
V$ are open {\bf E} and let $\Delta F = DF\Delta DF \hbox{\rm :}\
U \rightarrow L({\bf E}\Delta {\bf E}, {\bf E}\Delta {\bf E})$.
Then $D(\Delta F)\hbox{\rm :}\  U\rightarrow L({\bf E}, L({\bf
E}\Delta {\bf E}, {\bf E}\Delta {\bf E}))$ is given by
\begin{align*}
D(\Delta F)(u)(x, y\Delta z) &= [DF(u)\cdot y \Delta D^{2}F(u)(x
\Delta z)]\\[.3pc]
&\quad\, + [D^{2}F(u)(x\Delta y)\Delta DF(u)\cdot z ].
\end{align*}
\end{lem}

\begin{propo}$\left.\right.$\vspace{.5pc}

\noindent There is a one-to-one correspondence between third-order
Hessian structures on $M$ and special second-order connections on
$T^{(2)}M$ such that if $K$ and $\widetilde{\n }$ corresponds to
each other, then for every chart $(U, \va)$ on $M$, $X, Y, Z \in
\sx(M)$ and $f \in {\cal F}(U)$ we have $ Kf(X,Y,Z)=
D^{3}f_{\varphi}(X_{\varphi}, Y_{\varphi}, Z_{\varphi}) +
D^{2}f_{\varphi}\cdot C_{\varphi}(X_{\varphi}, Y_{\varphi},
Z_{\varphi}) + Df_{\varphi}\cdot R_{\varphi}(X_{\varphi},
Y_{\varphi}, Z_{\varphi})$, where $R_{\va }$ and $C_{\va}$ are the
Christoffel symbols of $\widetilde{\n}$ in the chart $(U, \va)$.
\end{propo}

\begin{proof}
Suppose $K$ is a third-order Hessian structure on $M$. Fix a chart
$(U,\varphi)$ on $M$. Fix $m \in U$ and $X, Y, Z \in \sx(U)$. Let
$\varphi m = u $. For $f\in {\cal F}(U)$, define $ \Lambda f(X, Y,
Z)(m)= Kf(X,Y,Z)(m) - D^{3}f_{\varphi}(X_{\varphi}, Y_{\varphi},
Z_{\varphi})(u). $ Then (i) $\Lambda $ is real linear in $f$ and
(ii) $\Lambda f(X,Y,Z)(m)=0$ if $Df(m) = 0$ and $D^{2}f(m)=0$. It
follows that $\Lambda \hbox{\rm :}\  f \rightarrow \Lambda
f(X,Y,Z)(m)$ is a second-order tangent vector at $m$. Hence there
exists a second-order tangent vector in ${\bf E}\oplus {\bf
E}\Delta {\bf E}$ which we  denote by $ R_{\varphi}(X_{\varphi},
Y_{\varphi}, Z_{\varphi})(u)\oplus C_{\varphi}(X_{\varphi},
Y_{\varphi}, Z_{\varphi})(u) $ such that $ \Lambda (f)=
 D^{2}f_{\varphi}\cdot C_{\varphi}(X_{\varphi}, Y_{\varphi}, Z_{\varphi})(u)
 + Df_{\varphi}\cdot R_{\varphi}(X_{\varphi}, Y_{\varphi}, Z_{\varphi})(u).
$ Therefore, if $X, Y, Z \in \sx(U)$, $ [R_{\varphi}(X_{\varphi},
Y_{\varphi}, Z_{\varphi})\oplus C_{\varphi}(X_{\varphi},
Y_{\varphi}, Z_{\varphi})](f_{\varphi}) = [Kf(X,Y,Z) \circ
\varphi^{-1} - D^{3}f_{\varphi}(X_{\varphi}, Y_{\varphi},
Z_{\varphi})]. $ We may conclude from this that $R_{\varphi}\oplus
C_{\varphi}$ is a smooth trilinear map on $\varphi U$. Therefore
$Kf(X,Y,Z)(m)= D^{3}f_{\varphi}(X_{\varphi}, Y_{\varphi},
Z_{\varphi})(u) + D^{2}f_{\varphi}\cdot C_{\varphi}(X_{\varphi},
Y_{\varphi}, Z_{\varphi})(u) + Df_{\varphi}\cdot
R_{\varphi}(X_{\varphi}, Y_{\varphi}, Z_{\varphi})(u).$ We can
easily verify the transformation laws for $C_{\varphi}$ and
$R_{\varphi}$ and hence we obtain $ C_{\psi}(X_{\psi}, Y_{\psi},
Z_{\psi}) = [DF\Delta DF]\cdot C_{\varphi}(X_{\varphi},
Y_{\varphi}, Z_{\varphi})(u) + [X_{\psi}(v)\Delta
D^{2}F(Y_{\varphi}, Z_{\varphi})(u) + \hbox{cyclic terms}]$ and $
R_{\psi}(X_{\psi}, Y_{\psi}, Z_{\psi}) = DF(u) \cdot
R_{\varphi}(X_{\varphi}, Y_{\varphi}, Z_{\varphi})(u) +
D^{2}F(u)\cdot C_{\varphi}(X_{\varphi}, Y_{\varphi},
Z_{\varphi})(u) + D^{3}F(X_{\varphi}, Y_{\varphi},
Z_{\varphi})(u).$ These are precisely the relations satisfied by
the Christoffel symbols ($C$ and $R$) of a special second-order
connection $\widetilde{\nabla}$, so that $K$ determines a special
second-order connection $\widetilde{\nabla}$ on $M$.

Conversely, if $\widetilde{\nabla}$ is a second-order connection
on $M$ and $C, R$ are the associated Christoffel symbols of
$\widetilde{\nabla}$ in a chart, then defining $Kf(X,Y,Z)(m)=
D^{3}f_{\varphi}(X_{\varphi}, Y_{\varphi}, Z_{\varphi})(u) +
D^{2}f_{\varphi}\cdot C_{\varphi}(X_{\varphi}, Y_{\varphi},
Z_{\varphi})(u) + Df_{\varphi}\cdot R_{\varphi}(X_{\varphi},
Y_{\varphi}, Z_{\varphi})(u)$, it is not difficult to verify, by
reversing the order of steps in the above proof, that $K$ defines
a third-order Hessian structure on $M$. Thus we have a one-to-one
correspondence between special second-order connections and
third-order Hessian structures.\hfill {\rm q.e.d.}
\end{proof}

In the next section we discuss the concept of a {\it geodesic} of
special second-order connection and we prove a theorem below which
gives a very satisfying geometric characterization of symmetric
third-order Hessian structures.\vspace{-.4pc}

\section{Second-order geodesics}

Let $c$ be a smooth curve in $M$ such that $c(0) = m$. For $f \in
{\cal F}(M)$, define $\ddot{c} (f) = (f \circ c)''(0)$. That is in
a local chart $(U, \va)$ at $m \in U$, we have $\ddot{c} (f) =
D^{2}f_{\va}(\va m) (c_{\va}'(0), c_{\va}'(0)) + Df_{\va}(\va
m)\cdot c_{\va}''(0)$. If $f\in F_{m}^{3}$ then $\ddot{c} (f) =
0$. Hence $\ddot{c} \in T_{m}^{(2)}M$. Let $(U, \va)$ be chart at
$m$. Define $[\ddot{c}]_{\va} \ \substack{{\rm def}\\ {=}} \ c_{\va}''
\oplus c_{\va}' \Delta c_{\va}'$. This is well-defined and is a
second-order tangent vector at $m$.

Note that, if $\widetilde{\n }$ is a special second-order
connection on $T^{(2)}M$ then $\widetilde{\n}_{\bf v}\tau $
depends only on the behavior of $\tau $ near $\pi ({\bf v})$ and
also that $\widetilde{\n}_{\bf v}\tau $ can be calculated once
$\tau $ is known along any smooth curve $c$ in $M$ with initial
tangent vector {\bf v} (see \cite{16}). Therefore
$\widetilde{\nabla}_{\dot{c}}\ddot{c}$ is well-\break defined.

Suppose $\widetilde{\nabla}$ is a special second-order connection
with Christoffel symbols $R$ and $C$. Let $(U, \varphi)$ be a
chart at $m$. Then we obtain
\begin{equation*}
(\widetilde{\nabla}_{\dot{c}}\ddot{c})_{\varphi} =
[c'''_{\varphi} - R_{\varphi}(c'_{\varphi}, c'_{\varphi},
c'_{\varphi})] \oplus [ C_{\varphi}(c'_{\varphi},
c'_{\varphi}, c'_{\varphi}) - 3 c'_{\varphi} \Delta
c''_{\varphi}].\vspace{.5pc}
\end{equation*}

\setcounter{theore}{0}
\begin{definit}$\left.\right.$\vspace{.5pc}

\noindent {\rm $c$ is called a {\it second-order geodesic} of
$\widetilde{\nabla}$ if $\widetilde{\nabla}_{\dot{c}}\ddot{c} = 0
$ for all $t$. Thus $c$ is a {\it second-order geodesic} if and
only if it satisfies the following equations:
\newpage

\begin{enumerate}
\renewcommand\labelenumi{(\roman{enumi})}
\item $c'''_{\varphi} = R_{\varphi}(c'_{\varphi}, c'_{\varphi},
c'_{\varphi})$ and

\item $C_{\varphi}(c'_{\varphi}, c'_{\varphi},
c'_{\varphi}) = 3 c'_{\varphi} \Delta c''_{\varphi}.$
\end{enumerate}}
\end{definit}

\begin{definit}$\left.\right.$\vspace{.5pc}

\noindent {\rm $\widetilde{\nabla}$ is said to be symmetric if $K$
is symmetric.}
\end{definit}

\begin{theor}[\!]
Let $\widetilde{\nabla}$ be a special second-order symmetric
connection on $T^{(2)}M$ and let $K$ be the unique symmetric
third-order Hessian structure on $M$ associated to
$\widetilde{\n}$. Then $Kf({\bf v, v, v}) = (f\circ c)'''(0)$
for all $m \in M${\rm ,} ${\bf v}\in T_{m}M$  if there is a
second-order geodesic $c$ such that $c(0) = m${\rm ,} $\dot{c} (0)
= {\bf v}$.
\end{theor}

(The proof is easy by using the definition and hence left to the
reader.)

Given a connection $\n$ on $TM$, there are standard theorems on
how $\n $ induces a connection on tensor bundles associated to $M$
by covariant differentiation. But we are not aware of any
observation in the literature that $\n$ induces a connection on
$T^{(2)}M$. We can now show that such a thing is possible.

\begin{definit}$\left.\right.$\vspace{.5pc}

\noindent {\rm Let $\n $ be a connection on $TM$. Then $Kf = \n
^{3} f$ defines a third-order Hessian structure on $M$ and
associated to $K$, there is a special second-order connection
$\widetilde{\n }$ on $T^{(2)}M$. We call $\widetilde{\n}$ the
connection induced by $\n$.} \end{definit}

\begin{propo}$\left.\right.$\vspace{.5pc}

\noindent Suppose $\n$ has the Christoffel symbol $\G_{\va}$ in a
chart $(U, \va )$. Then the induced connection $\widetilde{\n }$
on $T^{(2)}M$ has Christoffel symbols given by
\begin{align*}
C_{\va}(X_{\va}, Y_{\va}, Z_{\va}) &= X_{\va}\Delta
\G_{\va}(Y_{\va}, Z_{\va})
       + Y_{\va} \Delta \G_{\va}(X_{\va}, Z_{\va})
       + Z_{\va} \Delta \G_{\va}(X_{\va}, Y_{\va})\\
       \intertext{and}
R_{\va}(X_{\va}, Y_{\va}, Z_{\va}) &=  D[\G_{\va}(Y_{\va},
Z_{\va})]\. X_{\va}
       - \G_{\va}((\n_{X}Y)_{\va}, Z_{\va})\\
       &\quad\,-  \G_{\va}(Y_{\va}, (\n_{X}Z)_{\va}).
\end{align*}
\end{propo}

\begin{proof}
We know that in a local chart $(U, \va)$ at $m\in U$, $\n $ has
the form $(\n _{X}Y)_{\va} = DY_{\va}\cdot X_{\va} - \G
_{\va}(X_{\va}, Y_{\va}) $, where $\G_{\va}$ is the Christoffel
symbol of $\n $ which has transformation property \cite{5}. Then
in a local chart $(U, \va )$ at $m$, $\n ^{3}f$ has the expression
(and by using Leibnitz rule for higher derivatives, \cite{2})
\begin{align*}
\n ^{3}f(X, Y, Z) &= [\n _{X}(\n^{2}f)](Y, Z)\\[.2pc]
&= [\n _{X}(\n ^{2}f(Y, Z))] - [\n ^{2}f(\n _{X}Y, Z)]
- [ \n ^{2}f(Y, \n _{X}Z) ]\\[.2pc]
&=\l [ X\. (Y\. Z\. f - \n _{Y}Z\. f)\r ]
    - [ \n _{X}Y\. Z\. f - \n _{\n _{X}Y}Z\. f]\\[.2pc]
&\quad\, - \l [ Y\. \n _{X}Z\. f - \n_{Y}\n_{X}Z\. f \r ]\\[.2pc]
&= D^{3}f_{\va }(X_{\va}, Y_{\va}, Z_{\va}) + [D^{2}f_{\va}
(X_{\va}, \G_{\va}(Y_{\va}, Z_{\va}))\\[.2pc]
&\quad\, + D^{2}f_{\va} (Y_{\va}, \G_{\va}(X_{\va}, Z_{\va}))
       + D^{2}f_{\va} (Z_{\va}, \G_{\va}(X_{\va}, Y_{\va}))]\\[.2pc]
&\quad\, + [Df_{\va}\. D[\G_{\va}(Y_{\va}, Z_{\va})]\. X_{\va}
       - Df_{\va}\. \G_{\va}((\n_{X}Y)_{\va}, Z_{\va})\\[.2pc]
&\quad\, - Df_{\va}\.  \G_{\va}(Y_{\va}, (\n_{X}Z)_{\va})].
\end{align*}
That is,
\begin{align*}
\n ^{3}f(X, Y, Z) &= D^{3}f_{\va }(X_{\va}, Y_{\va}, Z_{\va})+
D^{2}f_{\va}\. [ X_{\va}\Delta \G_{\va}(Y_{\va}, Z_{\va})\\[.3pc]
&\quad\, + Y_{\va} \Delta \G_{\va}(X_{\va}, Z_{\va}) + Z_{\va}
\Delta \G_{\va}(X_{\va}, Y_{\va})]\\[.3pc]
&\quad\, + Df_{\va}\. [D[\G_{\va}(Y_{\va}, Z_{\va})]\. X_{\va}-
\G_{\va}((\n_{X}Y)_{\va}, Z_{\va})\\[.3pc]
&\quad\, - \G_{\va}(Y_{\va}, (\n_{X}Z)_{\va})].
\end{align*}

Hence we have $C_{\va}(X_{\va}, Y_{\va}, Z_{\va})$ and
$R_{\va}(X_{\va}, Y_{\va}, Z_{\va})$ as mentioned in Proposition
5.5.

\hfill {\rm q.e.d.}
\end{proof}

\begin{note}
{\rm If $\n $ has torsion and curvature zero, then $\n ^{3}f$ is
symmetric and hence the corresponding $\widetilde{\n }$ is
symmetric.}\vspace{-.4pc}
\end{note}

\section{Conclusion}

Not every $K$ arises like $\n ^{3}$, for some $\n$, whereas in
the second-order case, every $H^{2}$ arises like $\n^{2}$ for a
unique $\n$.

\setcounter{theore}{0}
\begin{propo}$\left.\right.$\vspace{.5pc}

\noindent If \ $\widetilde{\n}$ is a special second-order connection
on $T^{(2)}M$ induced from a (first order) connection $\n $ on
$M${\rm ,} then every  geodesic of $\n$ is also a geodesic of
$\widetilde{\n}$.
\end{propo}

\begin{proof}
Let $\n $ be a connection on $M$. Recall that $c$ is a geodesic of
$\n $ iff $c_{\va}'' = \G_{\va}(c_{\va}', c_{\va}')$, for all
charts $\va $. We have to prove that $c$ is a second-order
geodesic of $\widetilde{\n }$. Using the equations for $C_{\va}$
and $R_{\va}$ in Proposition 5.5, we obtain $C_{\va}(c_{\va}',
c_{\va}', c_{\va}') = 3 c_{\va}' \Delta \G_{\va}(c_{\va}',
c_{\va}')= 3 c_{\va}' \Delta c_{\va}''.$ and $R_{\va}(c_{\va}',
c_{\va}', c_{\va}') = c_{\va}'''.$ \hfill {\rm q.e.d.}\vspace{-.4pc}
\end{proof}

\section*{Acknowledgment}

I am deeply grateful to my guru Prof. K~Viswanath, Department of
Mathematics, University of Hyderabad, for suggesting this problem
as part of my Ph.D. thesis and my gratitude to him for his
continuous support and guidance.
\end{document}